\documentclass[11pt]{article}

\usepackage[a4paper, hmargin=25mm, top=28mm, bottom=30mm]{geometry}

\usepackage{amssymb,amsmath,amsthm}
\usepackage[T1]{fontenc}
\usepackage{lmodern}
\usepackage{microtype}

\usepackage[dvipsnames]{xcolor}
\usepackage{tikz-cd}
\usepackage[colorlinks,linkcolor=MidnightBlue,citecolor=OliveGreen,
  urlcolor=MidnightBlue]{hyperref}

\DeclareMathOperator{\Gal}{Gal}
\DeclareMathOperator{\Ind}{Ind}
\DeclareMathOperator{\Irr}{Irr}
\DeclareMathOperator{\Stab}{Stab}
\DeclareMathOperator{\cons}{cons}
\DeclareMathOperator{\Res}{Res}
\DeclareMathOperator{\core}{core}
\DeclareMathOperator{\GL}{GL}
\DeclareMathOperator{\PSL}{PSL}

\newcommand{\F}{\mathbb F}
\newcommand{\Q}{\mathbb Q}
\newcommand{\Z}{\mathbb Z}
\newcommand{\C}{\mathbb C}
\newcommand{\Zhat}{\widehat{\mathbb Z}}
\newcommand{\Gt}{\widetilde{G}}
\newcommand{\Ht}{\widetilde{H}}

\newtheorem{theorem}{Theorem}
\newtheorem{cor}[theorem]{Corollary}
\newtheorem{ex}[theorem]{Example}
\newtheorem{lem}[theorem]{Lemma}
\newtheorem{defin}[theorem]{Definition}
\newtheorem{prop}[theorem]{Proposition}
\newtheorem{rem}[theorem]{Remark}
\newtheorem*{theorem*}{Theorem}
\newtheorem*{question*}{Question}
\newcommand{\overbar}[1]{\mkern 1.5mu\overline{\mkern-1.5mu#1\mkern-1.5mu}\mkern 1.5mu}

\begin{document}
\title{A Note on Representations of Finite Groups with a View
  Towards the Neukirch--Uchida Theorem}
\author{Pavel Solomatin}
\date{2026}
\maketitle

\begin{abstract}
We prove two theorems on monomial representations of finite groups generalising results
of de Smit. Both give conditions under which a subgroup $H$ of a finite
group $G$ is determined up to conjugacy: the first by a single induced representation,
the second by the permutation representations of $H$ and of the kernel of a character
of $H$, or by the latter alone when this character has prime order. We then build a framework for theorems of Neukirch--Uchida
type on the second result and on a lemma on regular layers, extending an approach of an earlier paper of the author. A
profinite group whose isomorphic open subgroups have isomorphic permutation
representations has all isomorphic open subgroups conjugate if it has enough Clifford
characters, and every isomorphism between open subgroups is induced by conjugation if it
has enough regular layers. Applied to the absolute Galois group of $\Q$, this gives new
proofs of the Neukirch--Uchida theorem in its weak and strong forms.
\end{abstract}

\vfill\noindent\textbf{Use of AI tools.}
The author used Claude Opus 5.5 (Anthropic) for proofreading, editing, literature search
and computational checks in GAP and PARI/GP, and as a discussion partner. It also
contributed to the mathematics in several places, suggesting arguments, lemmas and
simplifications of proofs. The author has checked the arguments and takes
responsibility for the content.

\clearpage
\tableofcontents
\clearpage
\section{Introduction}\label{SecIntro}

\subsection{History}\label{SecIntroHistory}

Let $K$ be a number field with algebraic closure $\overbar K$. We write
$\mathcal G_K = \Gal(\overbar K:K)$ for its absolute Galois group, a profinite group,
and $\zeta_K$ for the Dedekind zeta function of $K$. Neukirch and Uchida
proved that a number field is determined by its absolute Galois group, and that every
isomorphism of absolute Galois groups comes from the fields:

\begin{theorem}[Neukirch, Uchida]\label{NUthm}
Let $K$ and $K'$ be number fields.
\begin{enumerate}
\renewcommand{\labelenumi}{(\alph{enumi})}
\item \emph{Weak form.} If $\mathcal G_K \simeq \mathcal G_{K'}$ as topological groups, then
$K \simeq K'$.
\item \emph{Strong form.} Let $K, K' \subseteq \overbar\Q$. Every isomorphism
$\varphi : \mathcal G_K \to \mathcal G_{K'}$ of topological groups is induced by
conjugation: there is $\sigma \in \mathcal G_\Q$ with $\sigma(K) = K'$ and
$\varphi(x) = \sigma x \sigma^{-1}$ for all $x \in \mathcal G_K$.
\end{enumerate}
\end{theorem}

The first version of the theorem was proved by Neukirch~\cite{NeukInvent} in 1969 under
the assumption that one of the two fields is Galois over $\Q$. His proof has two halves. The first
half recovers from $\mathcal G_K$ the decomposition groups of primes and with them the
splitting type of every rational prime $p$ in $K$, that is, the multiset of residue
degrees of the primes of $K$ above $p$. The second half recovers $K$ from these
splitting types when $K$ is Galois over $\Q$, which is essentially a theorem of Bauer~\cite[Ch.~VII, \S 13]{NeukANT}.

Uchida~\cite{Uchida1} removed the assumption that one field is Galois over $\Q$ by
proving the strong form directly and deducing the weak form. Uchida's
method extends to global function fields~\cite{Uchida2}. A modern account following Neukirch
and Uchida is in~\cite[Ch.~XII, \S 2]{neukirchCohomology}.
The theorem is the starting point of
anabelian geometry. Already the maximal prosolvable quotient of $\mathcal G_K$ suffices
to determine $K$ \cite{NeukAnotherProof},~\cite{Uchida3}.

A parallel story asks when number fields with equal splitting types are isomorphic. Let
$K$, $K'$ be number fields, $F$ a finite Galois extension of $\Q$ containing both with
$G = \Gal(F:\Q)$, and $H$, $H'$ the subgroups of $G$ fixing $K$ and $K'$. The zeta
function records the splitting types of the rational primes, and by Gassmann's criterion~\cite{Gassmann},~\cite{Perl} $\zeta_K = \zeta_{K'}$ if and only if the
permutation representations $\Ind_H^G(1) = \C[G/H]$ and $\Ind_{H'}^G(1) = \C[G/H']$
are isomorphic. Subgroups with isomorphic permutation representations are called
Gassmann equivalent. The smallest Gassmann equivalent subgroups which are not conjugate
have index $7$ \cite{Perl},~\cite{BosmadeSmit}. An example is
given by the two conjugacy classes of subgroups of $\PSL(2,7) \simeq \GL_3(\F_2)$
isomorphic to $S_4$. As $\PSL(2,7)$ is
a Galois group over $\Q$, they lead to non-isomorphic number fields with the same
splitting types.

This is where Neukirch meets Gassmann: the first half of his proof turns isomorphic Galois groups into
isomorphic permutation representations, and what remains of Theorem~\ref{NUthm}(a) is to
decide when subgroups with isomorphic permutation representations are conjugate.
For the arithmetic applications it is therefore natural to ask:
\begin{question*}
Which representations of a finite group $G$ determine a subgroup $H$ of $G$ up to
conjugacy?
\end{question*}

\subsection{Monomial representations}\label{SecIntroMonomial}

Looking for an invariant stronger than the zeta function, de Smit suggested monomial
representations, that is, representations induced from a character. He proved that
\emph{the isomorphism class of $K$ is determined by the $L$-function of a single
suitable character of $\mathcal G_K$}. By inductivity of Artin $L$-functions the
$L$-function $L(\chi,s)$ of such a character $\chi$ over $K$ is the $L$-function of
$\Ind_{\mathcal G_K}^{\mathcal G_\Q}(\chi)$ over $\Q$, and Artin $L$-functions over
$\Q$ determine the representation. So de Smit's statement is that this induced
representation determines $\mathcal G_K$ up to conjugacy. His argument works for
characters of any order $m \ge 3$ \cite[Thm.~10.1]{GlobalFieldsFromLFunct}. The
character comes from a wreath product of $C_m$ with the Galois group of the Galois
closure of $K$ (see Section~\ref{SecDeSmit}). Variations of the proof with explicit
trace computations appear in~\cite{SolNeuk},~\cite[Thm.~7]{ArEq} and~\cite[Ch.~3]{SoloThesis}, and in~\cite[Prop.~4.2.3]{CornelissenPeyerimhoff}, where
they are applied to isospectral Riemannian coverings.
A similar result is in the master's thesis of Pintonello~\cite{Pintonello},
supervised by de Smit: \cite[Thm.~2.2.2]{Pintonello} uses $C_2 \times C_2$ in place of
$C_m$. One quadratic $L$-function does not suffice for some fields, for instance
$\Q(\sqrt[8]{5})$ \cite[Thm.~3.2.2]{Pintonello}, while a suitable quadratic character together with $\zeta_K$ does~\cite[Thm.~2.3.1]{Pintonello}; see also~\cite[Rem.~4.3.3]{CornelissenPeyerimhoff}.

In this paper we generalise both group-theoretic approaches. Let $H$ be a subgroup of a finite group $G$ and $\chi$
a character of $H$. Theorem~\ref{General} says that the induced representation
$\Ind_H^G(\chi)$ determines the pair $(H,\chi)$ up to conjugacy when a certain orbit of
characters is line-free in the sense of Definition~\ref{DefLineFree}. A special case
is the group-theoretic part of de Smit's theorem (see Corollary~\ref{DeSmit}).

Theorem~\ref{CorFinite} uses permutation representations instead. It rests on the
stabiliser condition
\begin{equation*}
 \Stab_G(\chi|_{\core_G(H)}) \ = \ H , \tag{$\ast$}
\end{equation*}
where $\core_G(H)$ is the largest normal subgroup of $G$ contained in $H$ and
$\Stab_G$ is the stabiliser for the action of $G$ by conjugation on the characters of
$\core_G(H)$. We call a character $\chi$ satisfying $(\ast)$ a \emph{Clifford character}
of $H$, after the Clifford correspondence on which the proof rests (see Definition~\ref{DefClifford}). Put $U = \ker\chi$, and let $H'$ be a subgroup of $G$ and $U'$ a
normal subgroup of $H'$. Theorem~\ref{CorFinite} has two parts. Part
(i) says that if $H'/U'$ is abelian and the permutation representations of $H$ and
$U$ are isomorphic to those of $H'$ and $U'$, then $(H',U')$ is conjugate to $(H,U)$. Part (ii)
says that if $\chi$ has prime order $l$ and $[H':U'] = l$, then the isomorphism for $U$
and $U'$ alone suffices. Part (i) extracts what does not depend on the wreath product
from the arguments of~\cite[Thm.~3, Thm.~5]{SolNeuk} for $m = 3$, of~\cite[Thm.~2.1]{Phagan} for every $m \ge 3$ and of~\cite[Thm.~2.3.1]{Pintonello}
for $m = 2$. All of these use both permutation representations. The earlier proofs compute traces
explicitly and rely on the uniqueness of monomial structures. Our proofs use Mackey's
formula, Frobenius reciprocity and the Clifford correspondence~\cite[Thm.~6.11]{Isaacs}, and part (ii) also uses Galois conjugation of traces.

\subsection{A general Neukirch--Uchida theorem}\label{SecIntroGeneral}

In~\cite{SolNeuk} and~\cite[Ch.~2]{SoloThesis} a connection with Theorem~\ref{NUthm}(a) was proposed by replacing the $L$-function in de Smit's argument by the
zeta function of the cyclic extension of $K$ cut out by the character. An isomorphism
$\mathcal G_K \simeq \mathcal G_{K'}$ matches this extension with one of $K'$, and
the first half of Neukirch's proof, applied once more, shows that their zeta functions
agree. Phagan reworked the argument for every cyclic order $m \ge 3$
\cite[\S 2]{Phagan}. Sa\"{\i}di and Tamagawa prove~\cite[Thm.~1]{SaidiTamagawa} by
recovering decomposition groups and Frobenius elements and then applying~\cite{GlobalFieldsFromLFunct}, which is the same combination of the first half of
Neukirch's proof with de Smit's theorem. A forerunner of this approach is Uchida's proof of the
prosolvable case~\cite{Uchida3}, which already combines arithmetic equivalence of the
abelian extensions cut out by characters with the regular module. Using Theorem~\ref{CorFinite} and a lemma on regular layers we turn this approach into a general
framework.
We say that a profinite group $G$
\begin{itemize}
\item is \emph{conjugacy-rigid} if its isomorphic open subgroups are conjugate
(Definition~\ref{DefUchida});
\item is \emph{strongly conjugacy-rigid} if every isomorphism between open subgroups
is conjugation by an element of $G$ (Definition~\ref{DefStrong});
\item is \emph{Gassmann-faithful} if its isomorphic open subgroups have isomorphic
permutation representations (Definition~\ref{DefNeukirch});
\item has \emph{enough Clifford characters} if every open subgroup has a
Clifford character (Definition~\ref{DefInertial});
\item has \emph{enough regular layers} if every open normal subgroup $N_1$ of $G$
contains open normal subgroups $M \subseteq N$ of $G$ such that $N/M$ is abelian and
$N/M \simeq \F_l[G/N]$ as $G/N$-modules for a prime $l \ge 5$ (Definition~\ref{DefRegular}).
\end{itemize}

The main theorem of this paper is the following.

\begin{theorem*}[Theorem~\ref{MainThm}]
Let $G$ be a Gassmann-faithful profinite group.
\begin{enumerate}
\renewcommand{\labelenumi}{(\roman{enumi})}
\item If $G$ has enough Clifford characters, then $G$ is conjugacy-rigid.
\item If $G$ has enough regular layers, then $G$ is strongly conjugacy-rigid.
\end{enumerate}
\end{theorem*}

The main theorem applies to $\mathcal G_\Q$. Gassmann-faithfulness is Neukirch's
theorem in group-theoretic form: by~\cite{NeukInvent} isomorphic absolute Galois groups
of number fields have equal zeta functions, so by Gassmann's criterion $\mathcal G_\Q$
is Gassmann-faithful (see Proposition~\ref{ExGQNeukirch}). For part (i) a quadratic
character of each open subgroup suffices. Its kernel is built from a prime that splits
completely as in a construction of Pintonello~\cite[Prop.~2.2.1]{Pintonello} (see
Lemma~\ref{LemQuad}), and this gives Theorem~\ref{NUthm}(a). Part (ii) follows the
structure of Uchida's argument~\cite{Uchida1}, with a lemma on a single regular layer
in place of his final step (see Lemma~\ref{LemLayer}). We construct regular layers of
$\mathcal G_\Q$ by Kummer theory with $l = 5$ from the same kind of prime (see
Lemma~\ref{LemRegQ}), and this gives Theorem~\ref{NUthm}(b). The earlier approach
of~\cite{SolNeuk} realised de Smit's wreath product over $\Q$ by the Grunwald--Wang
theorem~\cite[Prop.~9.1]{GlobalFieldsFromLFunct}. Neither of our constructions needs
it or \v{S}afarevi\v{c}'s theorem on embedding problems, so, given the first half of
Neukirch's proof, we obtain both parts of Theorem~\ref{NUthm} without them.

The rest of the paper is organised as follows. Section~\ref{SecGroup} proves Theorems~\ref{General} and \ref{CorFinite}. Section~\ref{SecNeukUch} proves the main theorem and
Section~\ref{SecExamples} gives non-examples. Section~\ref{SecApplications} applies the
main theorem to $\mathcal G_\Q$ and to the Galois group of the maximal prosolvable
extension of $\Q$.

\section{Monomial representations of finite groups} \label{SecGroup}

\subsection{Notation}\label{SecNotation}

\paragraph{Groups and representations.}
All groups in this section are finite. Representations are finite dimensional over
$\C$. For a subgroup $H \le G$ and $g \in G$ we write ${}^gH = gHg^{-1}$, and
$\core_G(H) = \bigcap_{g \in G} {}^gH$. A \emph{character} of a group is a
homomorphism of finite order to $\C^\times$. We use the same letter for a
representation and its space, and for a character and the one-dimensional
representation it defines. Our reference for representation theory is~\cite{Serre}.
For $H \le G$ we write $\Ind_H^G$ and $\Res_H$ for induction and restriction~\cite[\S 3.3, \S 7.1]{Serre}.

The group $G$ acts on the finite set $G/H$, and $\Ind_H^G(1) = \C[G/H]$ is the
associated \emph{permutation representation}. Its trace is the \emph{permutation
character}, and its kernel is $\core_G(H)$, the kernel of the action of $G$ on $G/H$.
Two subgroups $H$, $H'$ of $G$ are \emph{Gassmann equivalent} if
$\Ind_H^G(1) \simeq \Ind_{H'}^G(1)$, or in other words if they have the same
permutation character. Isomorphic representations have equal kernels, so Gassmann
equivalent subgroups have equal cores.

For a character $\chi$ of $H$ and $g \in G$ put ${}^g(H,\chi) = ({}^gH, {}^g\chi)$ with
${}^g\chi(h) = \chi(g^{-1}hg)$ for $h \in {}^gH$. Conjugate pairs induce isomorphic
representations~\cite[\S 7.1]{Serre}.

For a finite abelian group $A$ we write $\Irr(A)$ for its group of characters under
pointwise multiplication; $|\Irr(A)| = |A|$. We write $C_k$ for a cyclic group of
order $k$.

\paragraph{Characters of a normal subgroup.}
For a group $G$ with a normal subgroup $M$ we write:
\begin{itemize}
\item ${}^g\lambda$ for a character $\lambda$ of $M$, as for pairs above; as
  ${}^gM = M$, this gives an action of $G$ on the characters of $M$;
\item $G\cdot\lambda$ and $\Stab_G(\lambda)$ for the orbit of $\lambda$ and
  $\{g \in G : {}^g\lambda = \lambda\}$;
\item $\cons(V) \subseteq \Irr(M)$ for $M$ abelian and a representation $V$ of a
  subgroup containing $M$, the set of characters $\chi_M$ of $M$ occurring in
  $\Res_M V$.
\end{itemize}

Our main tool is the Clifford correspondence. Let $N \trianglelefteq G$, let $\theta$
be a character of $N$ and $T = \Stab_G(\theta)$. Then $V \mapsto \Ind_T^G(V)$ is a
bijection from the irreducible representations $V$ of $T$ whose restriction to $N$ contains
$\theta$ onto the irreducible representations of $G$ whose restriction to $N$ contains
$\theta$ \cite[Thm.~6.11]{Isaacs}.

\subsection{Line-free orbits and the theorem on induced representations}
\label{SecLineFree}

The hypothesis of Theorem~\ref{General} below uses the following definition.

\begin{defin}\label{DefLineFree}
Let $A$ be a finite abelian group. A subset $T \subseteq A$ is \emph{line-free} if it
contains no coset of a non-trivial subgroup of $A$.
\end{defin}
If $A$ is an elementary abelian $l$-group, regarded as a vector space over $\F_l$ and
written additively, the cosets of subgroups of order $l$ are the affine lines
$v + \F_l u$ with $u \ne 0$, which explains the name.

The following fact about stabilisers is used repeatedly below. Let
$M \trianglelefteq G$ be a normal subgroup, not necessarily abelian. If
$M \le H' \le G$ and $\chi'$ is a character of $H'$, then
\begin{equation}\label{StabRes}
 H' \ \le \ \Stab_G(\chi'|_M) ,
\end{equation}
since $\chi'(u^{-1}cu) = \chi'(c)$ for $u \in H'$ and $c \in M$, as $\C^\times$ is
abelian. In particular
$M \le \Stab_G(\lambda)$ for every character $\lambda$ of $M$.

The first main result is the following theorem.

\begin{theorem}\label{General}
Let $G$ be a finite group, $M \trianglelefteq G$ an abelian normal subgroup,
$\lambda \in \Irr(M)$, $H = \Stab_{G}(\lambda)$, and let $\chi$ be a character
of $H$ with $\chi|_M = \lambda$. Suppose the orbit $O = G\cdot\lambda$ is a line-free
subset of $\Irr(M)$. Then $\Ind_H^{G}(\chi)$ is irreducible, and for every subgroup
$H' \le G$ and every character $\chi'$ of $H'$,
$$ \Ind_{H'}^{G}(\chi') \ \simeq \ \Ind_H^{G}(\chi)
\qquad\Longrightarrow\qquad (H',\chi') \ = \ {}^z(H,\chi) \ \text{ for some }
z \in G . $$
\end{theorem}

Such a character $\chi$ need not exist. Lemma~\ref{DeSmitConfig} gives the main
example.

\subsection{Preliminaries for the proof of \texorpdfstring{Theorem~\ref{General}}{the theorem on induced representations}}

We first sketch the proof.

\begin{proof}[Sketch of the proof of Theorem~\ref{General}]
Restricted to $M$ the representation $\Ind_H^G(\chi)$ is the sum of the characters in
the orbit $O$, each appearing once. We compare this with the restriction of
$\Ind_{H'}^G(\chi')$ to $M$.
If $M$ were not contained in $H'$, the characters of $M$ in this second restriction
would fill whole cosets of non-trivial subgroups of $\Irr(M)$ (see Lemma~\ref{Engine}). A
line-free orbit contains no such coset. So $M \le H'$. Then $\lambda$ occurs in the
restriction of $\Ind_{H'}^G(\chi')$ to $M$, so a conjugate of $H'$ lies in the
stabiliser $H$ of $\lambda$. It equals $H$ as the two induced representations have the
same dimension. The Clifford correspondence then identifies $\chi'$ with the
corresponding conjugate of $\chi$.
\end{proof}

To make this precise we need the following lemma. It describes the restriction of
$\Ind_H^G(\chi)$ to $M$ for an arbitrary subgroup $H$, not necessarily a stabiliser.

\begin{lem}\label{Engine}
Let $M \trianglelefteq G$ be a normal subgroup, $H \le G$ any subgroup and $\chi$ a
character of $H$, and put $B = M \cap H$.
\begin{enumerate}
\renewcommand{\labelenumi}{(\roman{enumi})}
\item If $M$ is abelian, then $\cons\big(\Ind_H^G(\chi)\big)$ is a non-empty union of
  cosets in $\Irr(M)$ of subgroups of order $[M:B]$.
\item If $M \le H$, then
  $\Res_M\Ind_H^G(\chi) \simeq \bigoplus_{tH \in G/H} {}^{t}(\chi|_M)$.
\item If $M \le H$ and a character $\lambda$ of $M$ occurs in
  $\Res_M\Ind_H^G(\chi)$, then there is $g \in G$ with ${}^g(\chi|_M) = \lambda$ and
  ${}^gH \le \Stab_G(\lambda)$.
\end{enumerate}
\end{lem}

\begin{proof}
By Mackey's formula~\cite[\S 7.3, Prop.~22]{Serre}
$$ \Res_M\Ind_H^G(\chi) \ \simeq \ \bigoplus_{t} \Ind_{M \cap {}^tH}^{M}({}^t\chi), $$
where $t$ runs over representatives of the double cosets $MtH$ in $G$, and
$M \cap {}^tH = {}^tB$ as $M$ is normal.

\medskip\noindent
\textbf{Proof of (i).} As $M$ is abelian, the summand of $t$ is a sum of characters of $M$
\cite[\S 3.1, Thm.~9]{Serre}, and by Frobenius reciprocity~\cite[\S 7.2, Thm.~13]{Serre} these are the characters restricting to ${}^t\chi$ on
${}^tB$. Restriction $\Irr(M) \to \Irr({}^tB)$ is a homomorphism whose kernel is the group of
characters trivial on ${}^tB$. So these characters form a coset of this kernel, which
is isomorphic to $\Irr(M/{}^tB)$ and has order $[M:{}^tB] = [M:B]$. The coset is non-empty as the
summand is non-zero.

\medskip\noindent
\textbf{Proof of (ii).} If $M \le H$, then ${}^tB = M$ and $MtH = tH$, so the summands are indexed by
$G/H$, and the summand of $t$ is ${}^t(\chi|_M)$.

\medskip\noindent
\textbf{Proof of (iii).} By (ii) we have $\lambda = {}^g(\chi|_M)$ for some $g \in G$. By
\eqref{StabRes} $H \le \Stab_G(\chi|_M)$, so
${}^gH \le \Stab_G({}^g(\chi|_M)) = \Stab_G(\lambda)$.
\end{proof}

\subsection{Proof of \texorpdfstring{Theorem~\ref{General}}{the theorem on induced
representations}}\label{SecProofGeneral}

\begin{proof}[Proof of Theorem~\ref{General}]
Let $M$, $\lambda$, $H$, $\chi$ and $O$ be as in Theorem~\ref{General}, put
$\psi = \Ind_H^G(\chi)$, let $H' \le G$, and let $\chi'$ be a character of $H'$
with $\psi' = \Ind_{H'}^G(\chi') \simeq \psi$.

\medskip\noindent
\textbf{Step 1. The constituents of $\psi$ on $M$.} Since
$M \le \Stab_G(\lambda) = H$, Lemma~\ref{Engine}(ii), applied to $H$ and $\chi$, gives
$\Res_M\psi \simeq \bigoplus_{tH \in G/H} {}^t\lambda$, where the ${}^t\lambda$ are distinct as $H = \Stab_G(\lambda)$, so $\cons(\psi) = O$. As $\chi$
is a character of $H = \Stab_G(\lambda)$ with $\chi|_M = \lambda$, the Clifford
correspondence shows that $\psi$ is irreducible.

\medskip\noindent
\textbf{Step 2. $M$ lies in $H'$.} Since $\psi \simeq \psi'$ we have
$\cons(\psi') = \cons(\psi) = O$, which is line-free. Put $B = M \cap H'$. By Lemma~\ref{Engine}(i), applied to $H'$ and $\chi'$, the set
$\cons(\psi')$ is a non-empty union of cosets of subgroups of $\Irr(M)$ of order
$[M:B]$. As $O$ is line-free, these subgroups are trivial. So $[M:B] = 1$, that is,
$M \le H'$.

\medskip\noindent
\textbf{Step 3. Conjugating $(H',\chi')$ to $(H,\chi)$.} By Step 2 we have $M \le H'$, and $\lambda$ lies in
$O = \cons(\psi')$. So Lemma~\ref{Engine}(iii), applied to $H'$ and $\chi'$, gives
$g \in G$ with ${}^g(\chi'|_M) = \lambda$ and ${}^gH' \le \Stab_G(\lambda) = H$. As
$\psi \simeq \psi'$ we have $[G:H'] = \dim\psi' = \dim\psi = [G:H]$, so
${}^gH' = H$. Then ${}^g\chi'$ and $\chi$ are characters of $H$ restricting to
$\lambda$ on $M$, and both induce $\psi$, as conjugate pairs induce isomorphic
representations (see Section~\ref{SecNotation}). By the Clifford correspondence for $M$
and $\lambda$, induction from $H$ to $G$ is injective on such characters. Hence
${}^g\chi' = \chi$, that is, $(H',\chi') = {}^z(H,\chi)$ with $z = g^{-1}$.
\end{proof}

\subsection{The wreath product and the realisation of a distinguished character}\label{SecDeSmit}

Following de Smit we build, for a finite group $G_0$, a subgroup $H_0$ and an integer
$m$, a group $\Gt$ with a subgroup $\Ht$ and a distinguished character $\chi$ of $\Ht$
which for $m \ge 3$ satisfies the hypotheses of Theorem~\ref{General}.

\paragraph{The wreath product.}
Fix an integer $m \ge 2$ and realise $C_m$ as the group $\mu_m(\C) \subset \C^\times$
of $m$-th roots of unity. Let $G_0$ be a finite group and $H_0 \le G_0$ a subgroup of
index $n$. Put
$$ M \ = \ C_m^{\,n} \ = \ \{\, c = (c_x)_{x \in G_0/H_0} \ : \
   c_x \in C_m \,\} ,$$
an abelian group of order $m^n$ written multiplicatively. It has one coordinate for
each of the $n$ cosets of $H_0$, and $G_0$ permutes the coordinates as it permutes the cosets:
\begin{equation}\label{Action}
 (g \cdot c)_x \ = \ c_{g^{-1}x} \qquad (g \in G_0, \ x \in G_0/H_0).
\end{equation}
The wreath product $\Gt = C_m \wr_{G_0/H_0} G_0$ is the semidirect product
$M \rtimes G_0$ for this action. Its elements are the pairs $(c,g)$ with $c \in M$ and
$g \in G_0$, and
$$ (c,g)(c',g') \ = \ (c\cdot(g\cdot c'),\,gg') .$$
We identify $M$ with $\{(c,1)\}$ and $G_0$ with $\{(1,g)\}$. Then $M$ is normal in
$\Gt$ and $(1,g)^{-1}(c,1)(1,g) = (g^{-1}\cdot c,1)$. Put
$$ \Ht \ = \ M \rtimes H_0 \ = \ \{(c,h) \ : \ c \in M, \ h \in H_0\} ,$$
a subgroup of index $n$ in $\Gt$.

\paragraph{The characters.}
For $x \in G_0/H_0$ let $\lambda_x \in \Irr(M)$ be the projection to the coordinate
$x$, that is, $\lambda_x(c) = c_x$. The distinguished character is the projection to
the coordinate of the trivial coset, extended to $\Ht$:
\begin{equation*}
 \chi : \Ht \longrightarrow C_m, \qquad \chi(c,h) \ = \ c_{H_0} .
\end{equation*}
It is a homomorphism since $(h\cdot c')_{H_0} = c'_{h^{-1}H_0} = c'_{H_0}$ for
$h \in H_0$. Put $\lambda = \chi|_M = \lambda_{H_0}$. We call
$(\Gt, \Ht, M, \lambda, \chi)$ the \emph{de Smit configuration} attached to
$(G_0, H_0, m)$. Its five entries are the data $G$, $H$, $M$, $\lambda$, $\chi$ of
Theorem~\ref{General}. The following diagram shows these subgroups of $\Gt$, with
$U = \ker\chi$. The edge from $\Gt$ to $\Ht$ carries the index, the other labels are
the quotients, and $\chi$ and $\lambda$ give $\Ht/U \simeq C_m$ and
$M/\ker\lambda \simeq C_m$.
\[
\begin{tikzcd}[column sep=0.8em, row sep=1.8em]
 & \Gt \arrow[d, dash, "n"] \arrow[ldd, dash, bend right=50, "G_0"']
 \arrow[dddd, dash, bend left=55, looseness=1.35, "C_m^{\,n} \rtimes G_0"] & \\
 & \Ht \arrow[ld, dash, "H_0"'] \arrow[rd, dash, "C_m"] & \\
M \arrow[rd, dash, "C_m"'] & & U \arrow[ld, dash, "H_0"] \\
 & \ker\lambda \arrow[d, dash, "C_m^{\,n-1}"'] & \\
 & 1 &
\end{tikzcd}
\]

\paragraph{The idea of the proof.}
Write $\Irr(M)$ additively, with
$\lambda_x$ the standard basis vector $e_x$ of $(\Z/m\Z)^n$. Conjugating $\lambda$ gives
the projections to all other coordinates, so the orbit of $\lambda$ is the standard
basis. A line through $e_x$ and $e_y$ contains $2e_x - e_y$, which is not a basis
vector unless $2 = 0$ in $\Z/m\Z$. The next lemma makes this precise and so verifies the
hypotheses of Theorem~\ref{General}.

\begin{lem}\label{DeSmitConfig}
Let $(\Gt, \Ht, M, \lambda, \chi)$ be the de Smit configuration attached to
$(G_0, H_0, m)$. Then $\lambda$ and $\chi$ have
order $m$, the orbit of $\lambda$ under $\Gt$ is $O = \{\lambda_x\}_{x \in G_0/H_0}$,
and the stabiliser of $\lambda$ in $\Gt$ is $\Ht$. If $m \ge 3$, then $O$ is
line-free.
\end{lem}

\begin{proof}
As above, $\Irr(M) = (\Z/m\Z)^n$ with $\lambda_x = e_x$, so $\lambda = e_{H_0}$ has order $m$, and so has $\chi$, which extends $\lambda$ and takes values in $C_m$.

\medskip\noindent
\textbf{Step 1. The orbit and the stabiliser of $\lambda$.} For
$c \in M$ and $g \in G_0$ the relation $(1,g)^{-1}(c,1)(1,g) = (g^{-1}\cdot c,1)$ and
\eqref{Action} give ${}^{(1,g)}\lambda_x(c) = (g^{-1}\cdot c)_x = c_{gx}$, that is,
${}^{(1,g)}\lambda_x = \lambda_{gx}$. As $M$ acts trivially on $\Irr(M)$, as noted
before Theorem~\ref{General}, and $(c,g) = (c,1)(1,g)$, we get
${}^{(c,g)}\lambda = \lambda_{gH_0}$. As $G_0$ acts transitively on $G_0/H_0$, the
orbit of $\lambda$ is $O$, and its stabiliser is $M \rtimes H_0 = \Ht$.

\medskip\noindent
\textbf{Step 2. Line-freeness for $m \ge 3$.} Let $m \ge 3$ and suppose
that $O$ contains a coset of a non-trivial subgroup $P$ of $\Irr(M)$. This coset is
$e_x + P$ for some $x$. Choose $u \ne 0$ in $P$. As $e_x + u \in O$, we have
$u = e_y - e_x$ for some $y \ne x$, so $e_x - u = 2e_x - e_y$ lies in $O$. Its
coordinates at $x$ and $y$ are $2$ and $-1$, which are non-zero as $m \ge 3$, so it is
not a basis vector, a contradiction.
\end{proof}

The bound $m \ge 3$ is necessary for line-freeness. For $m = 2$ any two distinct points
$v$ and $v+u$ of $\Irr(M) \simeq \F_2^{\,n}$ form a coset of $\{0,u\}$, so $O$ is
not line-free when $n \ge 2$. The conclusion of de Smit's theorem can then fail as
well, as the following example shows.

\begin{ex}\label{ExD4}
Take $(G_0,H_0,m) = (C_2,1,2)$. Then $\Gt = C_2 \wr C_2 = (C_2 \times C_2) \rtimes C_2$
is dihedral of order $8$, where the top $C_2$ swaps the two factors, and
$\Ht = M = C_2 \times C_2$. The character $\chi = \lambda$ is the projection onto the
first factor. The swap sends it to the projection
onto the second factor, so $\Ht$ is its stabiliser. Each of the other two subgroups of
order $4$ is not conjugate to $\Ht$ but carries a character inducing the same irreducible representation
as $\chi$ \cite[Rem.~10.3]{GlobalFieldsFromLFunct},~\cite[\S 2.1]{Pintonello}.
\end{ex}

For $m \ge 3$ we obtain de Smit's theorem~\cite[Thm.~10.1]{GlobalFieldsFromLFunct}.

\begin{cor}[de Smit]\label{DeSmit}
Let $m \ge 3$ and let $(\Gt, \Ht, M, \lambda, \chi)$ be the de Smit configuration
attached to $(G_0, H_0, m)$. Then $\Ind_{\Ht}^{\Gt}(\chi)$ is irreducible. Let
$H' \le \Gt$ be a subgroup and $\chi'$ a character of $H'$. If
$\Ind_{H'}^{\Gt}(\chi') \simeq \Ind_{\Ht}^{\Gt}(\chi)$, then
$(H',\chi') = {}^z(\Ht,\chi)$ for some $z \in \Gt$.
\end{cor}

\begin{proof}
By Lemma~\ref{DeSmitConfig} the de Smit configuration satisfies the hypotheses of
Theorem~\ref{General}.
\end{proof}

\subsection{Pairs of permutation representations}\label{SecZeta}

Theorem~\ref{General} uses one induced representation and needs a line-free orbit. We
now drop line-freeness and use permutation representations instead, together with the
following condition on $\chi$.

\begin{defin}\label{DefClifford}
Let $G$ be a finite group, $H \le G$ and $D = \core_G(H)$. A character $\chi$ of $H$ is a
\emph{Clifford character} of $H$ if
\begin{equation*}\tag{$\ast$}\label{Star}
 \Stab_G(\chi|_D) \ = \ H .
\end{equation*}
\end{defin}

Here $G$ acts on the characters of the normal subgroup $D$. By \eqref{StabRes} $H$
fixes $\chi|_D$, so $gH \mapsto {}^g(\chi|_D)$ is a well-defined map from $G/H$ to the
characters of $D$, and \eqref{Star} says that it is injective. The name refers to the Clifford correspondence for
$D$ of Section~\ref{SecNotation}, by which $\Ind_H^G(\chi)$ is irreducible. For
example, the character $\chi$ of Example~\ref{ExD4} is a Clifford character, as
$D = \Ht$ and the stabiliser of $\chi$ is $\Ht$.

The second main result is the following theorem.

\begin{theorem}\label{CorFinite}
Let $G$ be a finite group, $H \le G$ and $\chi$ a Clifford character of $H$. Put
$D = \core_G(H)$, $\lambda = \chi|_D$ and $U = \ker\chi$, and let $H' \le G$ and
$U' \trianglelefteq H'$.
\begin{enumerate}
\renewcommand{\labelenumi}{(\roman{enumi})}
\item If $H'/U'$ is abelian, $\Ind_H^G(1) \simeq \Ind_{H'}^G(1)$ and
  $\Ind_U^G(1) \simeq \Ind_{U'}^G(1)$, then $(H',U')$ is conjugate to $(H,U)$ in $G$.
\item If $\chi$ has prime order $l$, $[H':U'] = l$ and
  $\Ind_U^G(1) \simeq \Ind_{U'}^G(1)$, then $(H',U')$ is conjugate to $(H,U)$ in $G$.
\end{enumerate}
\end{theorem}

\begin{proof}[Proof of Theorem~\ref{CorFinite}(i)]
\textbf{Step 1. $\Ind_H^G(\chi)$ inside $\Ind_{U'}^G(1)$.} Gassmann equivalent subgroups have equal cores (see Section~\ref{SecNotation}), so $D = \core_G(H') \le H'$, and $|H| = |H'|$ and $|U| = |U'|$ by
comparing dimensions. By \eqref{Star} and the Clifford correspondence for $D$ the
representation $\Ind_H^G(\chi)$ is irreducible. As $\chi$ factors through $H/U$, it
occurs in $\C[H/U] = \Ind_U^H(1)$, so by transitivity of induction $\Ind_H^G(\chi)$ is
a constituent of $\Ind_U^G(1) \simeq \Ind_{U'}^G(1)$.

\medskip\noindent
\textbf{Step 2. Conjugating $H'$ to $H$.} As $H'/U'$ is abelian,
$\Ind_{U'}^{H'}(1) = \C[H'/U']$ is the sum of the characters $\rho$ of $H'/U'$, so
$\Ind_{U'}^G(1) = \bigoplus_\rho \Ind_{H'}^G(\rho)$. Hence $\Ind_H^G(\chi)$ is a
constituent of $\Ind_{H'}^G(\rho)$ for some $\rho$. By Lemma~\ref{Engine}(ii) $\lambda$
occurs in $\Res_D\Ind_H^G(\chi)$, so also in $\Res_D\Ind_{H'}^G(\rho)$, and Lemma~\ref{Engine}(iii), applied to $H'$ and $\rho$ with $M = D$, gives $g \in G$ with
${}^g(\rho|_D) = \lambda$ and ${}^gH' \le \Stab_G(\lambda) = H$ by \eqref{Star}. As
$|{}^gH'| = |H|$, we get ${}^gH' = H$.

\medskip\noindent
\textbf{Step 3. Conjugating $U'$ to $U$.} The characters ${}^g\rho$ and $\chi$
of $H$ both restrict to $\lambda$ on $D$. By the Clifford correspondence for $D$ the
representation $\Ind_H^G({}^g\rho) \simeq \Ind_{H'}^G(\rho)$ is irreducible, so it is
isomorphic to its constituent $\Ind_H^G(\chi)$. By the Clifford correspondence for $D$
and $\lambda$, induction from $H$ to $G$ is injective on such characters, so
${}^g\rho = \chi$. As $\rho$ is trivial on $U'$, the character
$\chi$ is trivial on ${}^gU'$, so ${}^gU' \le U$, with equality as $|U'| = |U|$.
\end{proof}

\begin{proof}[Proof of Theorem~\ref{CorFinite}(ii)]
\textbf{Step 1. Irreducibility of $\Ind_H^G(\chi)$.} As $[H:U] = l = [H':U']$ and $|U| = |U'|$ by comparing
dimensions, we have $|H| = |H'|$, and we put $n = [G:H] = [G:H']$. By \eqref{Star} and
the Clifford correspondence for $D$ the representation $\Ind_H^G(\chi)$ is irreducible
of dimension $n$.

\medskip\noindent
\textbf{Step 2. Matching $\Ind_H^G(\chi)$ with a representation induced from $H'$.} As $H'/U'$ has prime order $l$, it is cyclic, so $U' = \ker\chi'$ for a character
$\chi'$ of $H'$ of order $l$. As for $H'$ in the proof of (i), transitivity of induction gives
$$ \Ind_U^{G}(1) = \bigoplus_{j=0}^{l-1}\Ind_{H}^{G}(\chi^j),
\qquad
\Ind_{U'}^{G}(1)
  = \bigoplus_{j=0}^{l-1}\Ind_{H'}^{G}(\chi'^j). $$
The irreducible representation $\Ind_H^G(\chi)$ is a constituent of
$\Ind_U^G(1) \simeq \Ind_{U'}^G(1)$, so of some $\Ind_{H'}^G(\chi'^{j_0})$. This also
has dimension $n$, so the two are isomorphic. Moreover $j_0 \ne 0$: by Frobenius
reciprocity $\Ind_H^G(\chi)$ does not contain the trivial representation, as
$\chi \ne 1$, while $\Ind_{H'}^G(1)$ does. Replacing $\chi'$ by $\chi'^{j_0}$, which
also has kernel $U'$ and order $l$, we may assume $j_0 = 1$.

\medskip\noindent
\textbf{Step 3. Equal permutation representations of $H$ and $H'$.} The traces of $\Ind_H^G(\chi)$ and $\Ind_{H'}^G(\chi')$ take values in
$\Q(\mu_l)$. For $0 < k < l$ applying to them the automorphism $\sigma_k$ of
$\Q(\mu_l)$ with $\sigma_k(\varepsilon) = \varepsilon^k$ for $\varepsilon \in \mu_l$ gives
$\Ind_H^G(\chi^k) \simeq \Ind_{H'}^G(\chi'^k)$. Summing over $k$ and subtracting from
$\Ind_U^G(1) \simeq \Ind_{U'}^G(1)$ gives $\Ind_H^G(1) \simeq \Ind_{H'}^G(1)$, as
representations are determined by their traces~\cite[\S 2.3]{Serre}. Now $H'/U'$ is
abelian, and (i) applies.
\end{proof}

The setting of Theorem~\ref{General} always gives a Clifford character, without
line-freeness. If $M \trianglelefteq G$, $\lambda$ is a character of $M$,
$H = \Stab_G(\lambda)$ and $\chi$ extends $\lambda$, then $M \le H$ by
\eqref{StabRes}, so $M \le \core_G(H) = D$ as $M$ is normal. Hence
$\Stab_G(\chi|_D) \le \Stab_G(\lambda) = H$ and \eqref{Star} holds. This applies to
the de Smit configuration for every $m \ge 2$, as $\Ht$ is the stabiliser of $\lambda$
by Lemma~\ref{DeSmitConfig} and $\chi$ extends $\lambda$. So Theorem~\ref{CorFinite}(i) applies for every $m$, and (ii) for every
prime $m$, including $m = 2$ where the hypothesis of Theorem~\ref{General} fails when $n \ge 2$, as in Example~\ref{ExD4}. The
conclusion recovers the pair $(\Ht, \ker\chi)$, and hence $\chi$ only up to a character with the same kernel.

The next example shows that (ii) needs $l$ prime.

\begin{ex}\label{ExQ8}
In the
quaternion group $Q_8$ let $H = \langle i \rangle$, let $\chi$ be a faithful character
of $H$ of order $4$, and put $U = U' = 1$ and $H' = \langle j \rangle$. As $H$ is
normal, $D = H$ and $\lambda = \chi$. Conjugation by $j$ sends $\chi$ to $\chi^3$, so
\eqref{Star} holds. Trivially $\Ind_U^{Q_8}(1) \simeq \Ind_{U'}^{Q_8}(1)$ and
$[H':U'] = 4$. But $H'$ is not conjugate to $H$ \cite[Ex.~6]{ArEq}.
\end{ex}

The following remark shows that condition \eqref{Star} cannot be dropped.

\begin{rem}\label{RemStarNeeded}
Let $H$ and $H'$ be Gassmann
equivalent subgroups of $G$ which are not conjugate, and let $\psi$ be a character of
$G$ of prime order $l$ which is non-trivial on $H$. Such data exist, for instance
$G = \PSL(2,7) \times C_2$ with $H$ and $H'$ the products with $C_2$ of two
non-conjugate subgroups isomorphic to $S_4$, and $\psi$ the non-trivial character of $G$ that is trivial on $\PSL(2,7)$. Put $\chi = \psi|_H$,
$\chi' = \psi|_{H'}$, $U = \ker\chi$ and $U' = \ker\chi'$. By Frobenius reciprocity $\psi$
occurs in $\Ind_H^G(1)$ if and only if $\chi = 1$, and likewise for $H'$. So $\chi'$ is
non-trivial and $[H':U'] = l$. The projection formula~\cite[\S 3.3]{Serre} gives
$$ \Ind_{H}^{G}(\chi^j)
 \simeq \psi^j \otimes \Ind_{H}^{G}(1)
  \simeq \psi^j \otimes \Ind_{H'}^{G}(1)
  \simeq \Ind_{H'}^{G}(\chi'^j) $$
for all $j$. Summing over $0 \le j < l$ as in the proof of Theorem~\ref{CorFinite}(ii) gives $\Ind_U^G(1) \simeq \Ind_{U'}^G(1)$. Thus all hypotheses of both parts of Theorem~\ref{CorFinite} hold except \eqref{Star}, which fails as $\chi|_D = \psi|_D$ is invariant under
$G$, where $D = \core_G(H)$. But $(H',U')$ is not conjugate to $(H,U)$, as $H'$ is not conjugate to $H$. This use of the projection formula goes
back to Nagata~\cite{Nagata}.
\end{rem}

Theorem~\ref{CorFinite} needs a large core. If $\core_G(H) = 1$ then $\lambda = 1$ has
stabiliser $G$, so \eqref{Star} fails unless $H = G$. Thus Theorem~\ref{CorFinite} never separates
core-free subgroups such as the Gassmann equivalent subgroups of $\PSL(2,7)$ within
$\PSL(2,7)$. This is why we pass to larger groups with large cores:
de Smit's wreath product (see Section~\ref{SecDeSmit}) and profinite groups such as
$\mathcal G_\Q$, where the core of $\mathcal G_K$ is $\mathcal G_N$ for the Galois
closure $N$ of $K$ (see the proof of Proposition~\ref{PropInertialQ}).

\section{Conjugacy-rigid groups and the main theorem} \label{SecNeukUch}

We now ask when isomorphic open subgroups of a profinite group must be conjugate.
From now on groups may be profinite. Homomorphisms between profinite groups are
continuous, and so are characters. Subgroups are closed. The subgroups we compare are
always open. A representation of a profinite group $G$ means one factoring through a
finite quotient. For an open normal subgroup $W$ of $G$, composing a representation of
$G/W$ with the projection $G \to G/W$ gives a representation of $G$, said to be \emph{inflated} from
$G/W$. Let $H$ be
open in $G$ and $\chi$ a character of $H$. Then $\ker\chi$ is open, so it contains an
open normal subgroup $W$ of $G$, for instance $W = \core_G(\ker\chi)$, which is open as
$\ker\chi$ has finitely many conjugates. We define
$\Ind_H^G(\chi)$ as the representation inflated from $\Ind_{H/W}^{G/W}(\chi)$; it does
not depend on $W$. In this way the notation and facts of Section~\ref{SecGroup} apply
to open subgroups of profinite groups.

\subsection{Conjugacy-rigid and Gassmann-faithful groups}

Theorem~\ref{NUthm}(a) follows once we know that $\mathcal G_\Q$ has the following
property (see Section~\ref{SecProofNU}).

\begin{defin}\label{DefUchida}
We call a profinite group $G$ \emph{conjugacy-rigid} if for every pair of open subgroups
$H$, $H' \le G$,
$$ H \simeq H' \ \text{ as profinite groups } \quad \Longrightarrow \quad
H \text{ and } H' \text{ are conjugate in } G. $$
\end{defin}

A stronger property asks the isomorphism itself to be a conjugation.

\begin{defin}\label{DefStrong}
We call a profinite group $G$ \emph{strongly conjugacy-rigid} if for every
isomorphism $\varphi : H \to H'$ of profinite groups between open subgroups of $G$ there is $\tau \in G$
with $\varphi(x) = \tau x \tau^{-1}$ for all $x \in H$.
\end{defin}

Then $H' = \varphi(H) = {}^\tau H$, so every strongly conjugacy-rigid group is conjugacy-rigid.
In the other direction, a weaker property asks isomorphic open subgroups only to be
Gassmann equivalent.

\begin{defin}\label{DefNeukirch}
We call a profinite group $G$ \emph{Gassmann-faithful} if for every pair of open
subgroups $H$, $H' \le G$,
$$ H \simeq H' \ \text{ as profinite groups } \quad \Longrightarrow \quad
\Ind_H^G(1) \simeq \Ind_{H'}^G(1). $$
\end{defin}

Conjugate subgroups are Gassmann equivalent, as conjugate pairs induce isomorphic
representations (see Section~\ref{SecNotation}). So every conjugacy-rigid group is
Gassmann-faithful.

\subsection{Enough Clifford characters}

For $H$ open in $G$ the core $\core_G(H)$ is open and normal in $G$, since $H$ has
finitely many conjugates, so Clifford characters of $H$ are defined as in Definition~\ref{DefClifford}. Part (i) of our main theorem uses the following property.

\begin{defin}\label{DefInertial}
A profinite group $G$ \emph{has enough Clifford characters} if every open
subgroup $H$ of $G$ has a Clifford character.
\end{defin}

Clifford characters pass to quotients by normal subgroups in their kernel, as the
following lemma shows.

\begin{lem}\label{LemCliffQuot}
Let $H$ be an open subgroup of a profinite group $G$, $\chi$ a character of $H$ and $W$ a closed normal
subgroup of $G$ contained in $\ker\chi$. Then $\chi$ is a Clifford character of $H$ if
and only if the character it defines on $H/W$ is a Clifford character of $H/W$ in
$G/W$.
\end{lem}

\begin{proof}
As $W$ is normal and lies in $\ker\chi \subseteq H$, it lies in $D = \core_G(H)$, and the core of $H/W$ in
$G/W$ is $D/W$. An element $g \in G$ fixes $\chi|_D$ if and only if $gW$ fixes the
character it defines on $D/W$. So the stabiliser of this character in $G/W$ is
$\Stab_G(\chi|_D)/W$, and it equals $H/W$ if and only if $\Stab_G(\chi|_D) = H$.
\end{proof}

\subsection{Enough regular layers}

Part (ii) of our main theorem replaces Clifford characters by regular layers. A
\emph{layer} of a profinite group $G$ is a quotient $N/M$ of open normal subgroups
$M \subseteq N$ of $G$. For a finite group $\Gamma$ and a prime $l$ the \emph{regular
module} $\F_l[\Gamma]$ has basis $(e_x)_{x \in \Gamma}$ and action
$g \cdot e_x = e_{gx}$. We call a layer $N/M$ \emph{regular} for a prime $l$ if
$N/M$ is abelian and $N/M \simeq \F_l[G/N]$ as $G/N$-modules, with $G/N$ acting by conjugation. The following definition asks for
such layers inside every open normal subgroup.

\begin{defin}\label{DefRegular}
A profinite group $G$ \emph{has enough regular layers} if for every open normal
subgroup $N_1$ of $G$ there are open normal subgroups $M \subseteq N \subseteq N_1$ of
$G$ such that $N/M$ is regular for a prime $l \ge 5$.
\end{defin}

The proof of part (ii) uses the following elementary fact on Gassmann equivalence.

\begin{lem}\label{LemGassCyclic}
Let $H$ and $H'$ be Gassmann equivalent subgroups of a finite group $G$. If $H$ is
cyclic, then $H'$ is conjugate to $H$.
\end{lem}

\begin{proof}
Comparing dimensions gives $|H'| = |H|$. Let $H = \langle h \rangle$. The permutation
character of $G$ on $G/H'$ takes at $h$ the number of cosets $xH'$ with
$x^{-1}hx \in H'$. It equals the permutation character on $G/H$ at $h$, which is
non-zero as $h$ fixes the coset $H$. So $x^{-1}hx \in H'$ for some $x \in G$, and as
$|H'| = |H|$ we get $H' = \langle x^{-1}hx \rangle = {}^{x^{-1}}H$.
\end{proof}

The heart of the argument is the following statement about a single layer.

\begin{lem}\label{LemLayer}
Let $G$ be a finite group and $M$ an abelian normal subgroup of $G$. Put
$\Gamma = G/M$ and suppose that $M \simeq \F_l[\Gamma]$ as $\Gamma$-modules for a
prime $l \ge 5$, with $\Gamma$ acting by conjugation. Let $H$ and $H'$ be subgroups of $G$ containing $M$ and let
$\psi : H \to H'$ be an isomorphism with $\psi(M) = M$. Suppose that $\psi(B)$ is
conjugate in $G$ to $B$ for every cyclic subgroup $B$ of $M$. Then there is
$h \in \Gamma$ with $\bar\psi(\delta) = h\delta h^{-1}$ for all $\delta \in H/M$,
where $\bar\psi : H/M \to H'/M$ is induced by $\psi$.
\end{lem}

\begin{proof}
\textbf{Step 1. The restriction of $\psi$ to $M$.} We write $M$ additively and identify it with $\F_l[\Gamma]$. The restriction $T$ of
$\psi$ to $M$ is additive and hence $\F_l$-linear. For $x \in H$ and $a \in M$ we have
$\psi(xax^{-1}) = \psi(x)\psi(a)\psi(x)^{-1}$, so $T(\delta \cdot a) =
\bar\psi(\delta) \cdot T(a)$ for $\delta \in H/M$. The conjugates of
$\langle w \rangle$ in $G$ are the subgroups $\langle g \cdot w \rangle$ with
$g \in \Gamma$. So for every $w \in M$ there are $g \in \Gamma$ and
$c \in \F_l^\times$ with $T(w) = c\, g \cdot w$.

\medskip\noindent
\textbf{Step 2. An element $h$ with $\bar\psi(\delta) = h\delta h^{-1}$ for all $\delta$.} For $w = e_1$ this gives $T(e_1) = c\, e_h$ with $c \in \F_l^\times$ and
$h \in \Gamma$. Hence $T(e_\delta) = T(\delta \cdot e_1) = c\, e_{\bar\psi(\delta)h}$
for $\delta \in H/M$. Now let $\delta \ne 1$ and put $w = e_1 + 2e_\delta$. Then
$$ c\,(e_h + 2e_{\bar\psi(\delta)h}) = c'\,(e_x + 2e_{x\delta}) $$
for some $x \in \Gamma$ and $c' \in \F_l^\times$. As $l \ne 2$ and
$\bar\psi(\delta) \ne 1$ both sides have support of size two. So either $h = x$ or
$h = x\delta$. In the second case $\bar\psi(\delta)h = x$, and comparing
coefficients gives $c = 2c'$ and $2c = c'$. Then $3c = 0$, which is impossible as
$l \ne 3$. Hence $h = x$ and $\bar\psi(\delta)h = h\delta$, that is,
$\bar\psi(\delta) = h\delta h^{-1}$.
\end{proof}

Lemma~\ref{LemLayer} fails for $l = 3$. Take $G = C_3 \wr C_3$ with
$M = \F_3[C_3]$ and the automorphism which inverts the top group $C_3$ and sends
$e_x$ to $e_{x^{-1}}$. With $H = H' = G$ this automorphism $\psi$ satisfies
$\psi(M) = M$ and maps every subgroup of $M$ to a conjugate. So it satisfies the
hypotheses of the lemma, but it is not inner modulo $M$. This does not show that Theorem~\ref{MainThm}(ii) would fail with $l = 3$ allowed in Definition~\ref{DefRegular}.

\subsection{The main theorem}

The following theorem separates the group theory of the Neukirch--Uchida theorem from
its arithmetic inputs. Part (i) is a formal consequence of Theorem~\ref{CorFinite}(i).

\begin{theorem}\label{MainThm}
Let $G$ be a Gassmann-faithful profinite group.
\begin{enumerate}
\renewcommand{\labelenumi}{(\roman{enumi})}
\item If $G$ has enough Clifford characters, then $G$ is conjugacy-rigid.
\item If $G$ has enough regular layers, then $G$ is strongly conjugacy-rigid.
\end{enumerate}
\end{theorem}

\begin{proof}[Proof of Theorem~\ref{MainThm}(i)]
Let $G$ be as in part (i) and let $\varphi : H \to H'$ be an isomorphism between
open subgroups of $G$. Let $\chi$ be a Clifford character of $H$, which exists by
Definition~\ref{DefInertial}. Put $U = \ker\chi$ and
$U' = \varphi(U)$. As $\varphi$ is a homeomorphism, $U'$ is open. It is normal in
$H'$, and $H'/U' \simeq H/U$ is abelian.
Since $H \simeq H'$ and $U \simeq U'$ and $G$ is Gassmann-faithful, we get
$$ \Ind_H^G(1) \simeq \Ind_{H'}^G(1), \qquad \Ind_U^G(1) \simeq \Ind_{U'}^G(1). $$
Let $W$ be an open normal subgroup of $G$ contained in $U \cap U'$, for instance $\core_G(U \cap U')$. The four
permutation representations are inflated from $G/W$, and $\chi$ defines a Clifford
character of $H/W$ in $G/W$ by Lemma~\ref{LemCliffQuot}. Theorem~\ref{CorFinite}(i) for $G/W$ gives $g \in G$ with $({}^gH)/W = H'/W$, so
${}^gH = H'$.
\end{proof}

This proof gives some $g$ with ${}^gH = H'$, but not one for which $\varphi$ is
conjugation by $g$. Part (ii) gives such an element. Its proof follows Uchida's proof~\cite{Uchida1}, whose final step needs a relation between $l$ and the order of the group acting on the layer, and \v{S}afarevi\v{c}'s theorem on embedding problems~\cite{Safarevic} to realise the layer. We replace this step by Lemma~\ref{LemLayer}, which works for every prime $l \ge 5$.

\begin{proof}[Proof of Theorem~\ref{MainThm}(ii)]
Let $\varphi : H \to H'$ be an isomorphism between open subgroups of $G$. As $\varphi$ is a homeomorphism, every open subgroup $V$ of $H$ is isomorphic to the open subgroup $\varphi(V)$ of $G$, so $V$ and $\varphi(V)$ are Gassmann equivalent as $G$ is Gassmann-faithful.

\medskip\noindent
\textbf{Step 1. A directed family of open normal subgroups.} Let $\mathcal F$ be the set of open normal subgroups $N \subseteq H$ of $G$ for which
there is an open normal subgroup $M \subseteq N$ of $G$ such that $N/M$ is regular for
a prime $l \ge 5$. Applying Definition~\ref{DefRegular} to
$N_1 \cap \core_G(H)$ shows that every open normal subgroup $N_1$ of $G$ contains a
member of $\mathcal F$. So $\mathcal F$ is directed under reverse inclusion and its
intersection is trivial.

\medskip\noindent
\textbf{Step 2. Conjugation modulo each $N \in \mathcal F$.} Fix $N \in \mathcal F$ with $M$ and $l$ as above. The subgroups $N$ and
$\varphi(N)$ are Gassmann equivalent, so they have the same index and the same core
(see Section~\ref{SecNotation}). As $N$ is normal, $N = \core_G(\varphi(N)) \subseteq
\varphi(N)$, and hence $\varphi(N) = N$. In the same way $\varphi(M) = M$. In
particular $N \subseteq H'$.
Then $\varphi$ induces an isomorphism $\psi : H/M \to H'/M$ with $\psi(N/M) = N/M$. Let
$B \le N/M$ be cyclic with preimage $V$ in $G$. Then $V$ and $\varphi(V)$ contain $M$,
so their permutation representations are inflated from $G/M$. Hence $B = V/M$ and
$\psi(B) = \varphi(V)/M$ are Gassmann equivalent in $G/M$, and they are conjugate by
Lemma~\ref{LemGassCyclic}. Let $\varphi_N : H/N \to H'/N$ be induced by $\varphi$
and put
$$ S_N = \{ h \in G/N : \varphi_N(\delta) = h\delta h^{-1} \text{ for all }
\delta \in H/N \} . $$
Lemma~\ref{LemLayer}, applied to the finite group $G/M$ with $N/M$ in place of $M$,
$\Gamma = G/N$ and $H/M$, $H'/M$ in place of $H$, $H'$, shows that the finite set $S_N$
is not empty, as $\varphi_N$ is the map $\bar\psi$ of the lemma.

\medskip\noindent
\textbf{Step 3. Conjugation by an element of $G$.} For $N' \subseteq N$ in $\mathcal F$ the projection $G/N' \to G/N$ maps
$S_{N'}$ into $S_N$. An inverse limit of non-empty finite sets over a directed
set is not empty. As every open normal subgroup of $G$ contains a member of $\mathcal F$, we have $G = \varprojlim_{N \in \mathcal F} G/N$, so there is $\tau \in G$ whose
image lies in $S_N$ for every $N \in \mathcal F$. For $x \in H$ this gives
$\varphi(x) \equiv \tau x \tau^{-1}$ modulo every $N \in \mathcal F$, and hence
$\varphi(x) = \tau x \tau^{-1}$.
\end{proof}

\section{Non-examples} \label{SecExamples}

This section shows that some natural profinite groups are not Gassmann-faithful. So
for them the first hypothesis of Theorem~\ref{MainThm} already fails, and as every
conjugacy-rigid group is Gassmann-faithful they are not conjugacy-rigid. Both examples rest on one observation. A profinite group $G$ with a
proper open subgroup isomorphic to $G$ is not Gassmann-faithful, as Gassmann
equivalent subgroups have equal index.

The simplest example is the absolute Galois group of a finite field.

\begin{ex}[$\Zhat = \mathcal G_{\F_q}$]\label{ExZhatFails}
For a finite field $\F_q$ the Frobenius $x \mapsto x^q$ generates
$\mathcal G_{\F_q}$ topologically, giving
$\mathcal G_{\F_q} \simeq \Zhat = \varprojlim_m \Z/m\Z$; the open subgroup of index
$m$ is $m\Zhat = \mathcal G_{\F_{q^m}}$. Multiplication by $m$ is an isomorphism
$\Zhat \to m\Zhat$. For $m \ge 2$ this subgroup is proper, so $\Zhat$ is not
Gassmann-faithful. In field terms, every
finite extension of $\F_q$ has absolute Galois group isomorphic to
$\mathcal G_{\F_q}$, so the analogue of Theorem~\ref{NUthm}(a) fails for finite fields.
The same argument applies to $\Z_p$, since $p\Z_p$ is a proper open subgroup isomorphic
to it.
\end{ex}

The maximal abelian quotient $\mathcal G_\Q^{\mathrm{ab}}$ of $\mathcal G_\Q$ fails in
the same way.

\begin{ex}[$\mathcal G_\Q^{\mathrm{ab}}$]\label{ExAbelian}
By the Kronecker--Weber theorem, the cyclotomic character identifies
$\mathcal G_\Q^{\mathrm{ab}}$ with $\Zhat^\times = \prod_p \Z_p^\times$. For an odd
prime $p$ the group $\Z_p^\times$ is the product of a cyclic group of order $p-1$ and
$1 + p\Z_p \simeq \Z_p$. Replacing the factor $1 + p\Z_p$ by its open subgroup $1 + p^2\Z_p \simeq p\Z_p$, as in
Example~\ref{ExZhatFails}, gives a proper open subgroup of $\mathcal G_\Q^{\mathrm{ab}}$
isomorphic to $\mathcal G_\Q^{\mathrm{ab}}$, so $\mathcal G_\Q^{\mathrm{ab}}$ is not
Gassmann-faithful.
\end{ex}

\section{Applications} \label{SecApplications}

\subsection{The Neukirch--Uchida theorem}\label{SecProofNU}

We first apply Theorem~\ref{MainThm} to $\mathcal G_\Q$. By Galois theory the open
subgroups of $\mathcal G_\Q$ are the groups $\mathcal G_K$ for number fields
$K \subseteq \overbar\Q$, and $[\mathcal G_\Q : \mathcal G_K] = [K:\Q]$. The first
hypothesis follows from a theorem of Neukirch.

\begin{prop}[Gassmann-faithful]\label{ExGQNeukirch}
$\mathcal G_\Q$ is Gassmann-faithful.
\end{prop}

\begin{proof}
Let $\mathcal G_K \simeq \mathcal G_{K'}$. By a theorem of Neukirch~\cite{NeukInvent},~\cite[Ch.~XII, \S 2]{neukirchCohomology} we have $\zeta_K = \zeta_{K'}$. So for a
finite Galois extension $F$ of $\Q$ containing $K$ and $K'$ Gassmann's criterion makes
$\Gal(F:K)$ and $\Gal(F:K')$ Gassmann equivalent in $\Gal(F:\Q)$, and
$\Ind_{\mathcal G_K}^{\mathcal G_\Q}(1)$ and $\Ind_{\mathcal G_{K'}}^{\mathcal G_\Q}(1)$
are inflated from their isomorphic permutation representations.
\end{proof}

For the second hypothesis a quadratic character suffices. Its kernel is
$\mathcal G_{K(\sqrt a)}$ for an element $a$ of $K$ as in the following lemma, a variant of a construction of Pintonello~\cite[Prop.~2.2.1]{Pintonello}.

\begin{lem}\label{LemQuad}
Let $K$ be a number field with Galois closure $N$. There is $a \in K^\times$, not a
square in $K$, such that for $g \in \Gal(N:\Q)$, $g(a)/a \in N^{\times 2}$ implies
$g \in \Gal(N:K)$.
\end{lem}

\begin{proof}
Let $p$ be a prime that splits completely in $N$, which exists by the Chebotarev
density theorem~\cite[Ch.~VII, Thm.~13.4]{NeukANT}. Fix a prime $\mathfrak P$ of $N$
above $p$ and put $\mathfrak p = \mathfrak P \cap \mathcal O_K$. By the Chinese
remainder theorem there is $a \in \mathcal O_K$ with $v_{\mathfrak p}(a) = 1$ and
$v_{\mathfrak q}(a) = 0$ for the other primes $\mathfrak q$ of $K$ above $p$. As
$v_{\mathfrak p}(a)$ is odd, $a$ is not a square in $K$. As $p$ is unramified in $N$, a
prime $g\mathfrak P$ with $g \in \Gal(N:\Q)$ has $v_{g\mathfrak P}(a) = 1$ if it lies
above $\mathfrak p$ and $v_{g\mathfrak P}(a) = 0$ otherwise. Now let $g(a)/a$ be a
square in $N$. Then
$$ v_{g\mathfrak P}(a) \ \equiv \ v_{g\mathfrak P}(g(a)) \ = \ v_{\mathfrak P}(a)
   \ = \ 1 \pmod 2 , $$
so $g\mathfrak P$ lies above $\mathfrak p$. The primes of $N$ above $\mathfrak p$ are
the $h\mathfrak P$ with $h \in \Gal(N:K)$, so $g\mathfrak P = h\mathfrak P$ for such
an $h$. As $p$ splits completely in $N$ the decomposition group of $\mathfrak P$ is
trivial, and hence $g = h \in \Gal(N:K)$.
\end{proof}

With it we show that $\mathcal G_\Q$ has enough Clifford characters.

\begin{prop}[Enough Clifford characters]\label{PropInertialQ}
$\mathcal G_\Q$ has enough Clifford characters.
\end{prop}

\begin{proof}
Let $H = \mathcal G_K$ be an open subgroup and $N$ the Galois closure of $K$. The
conjugates of $\mathcal G_K$ are the groups $\mathcal G_{g(K)}$ for
$g \in \mathcal G_\Q$, and $N$ is the compositum of the fields $g(K)$. So
$\core_{\mathcal G_\Q}(H) = \mathcal G_N$. Take $a$ as in Lemma~\ref{LemQuad} and let
$\chi$ be the character of $\mathcal G_K$ with kernel $\mathcal G_{K(\sqrt a)}$, which
has order $2$ as $a$ is not a square in $K$. Its restriction $\lambda$ to
$\mathcal G_N$ has kernel $\mathcal G_{N(\sqrt a)}$, so for $g \in \mathcal G_\Q$ the
kernel of ${}^g\lambda$ is $g\mathcal G_{N(\sqrt a)}g^{-1} = \mathcal G_{N(\sqrt{g(a)})}$.
Characters of order at most $2$ are determined by their kernels. So if
${}^g\lambda = \lambda$ then $N(\sqrt{g(a)}) = N(\sqrt a)$, and $g(a)/a \in N^{\times 2}$
by Kummer theory. By Lemma~\ref{LemQuad}, applied to the restriction of $g$ to $N$, this
implies $g \in \mathcal G_K$. So $\Stab_{\mathcal G_\Q}(\lambda) \subseteq \mathcal G_K$,
and the reverse inclusion is \eqref{StabRes}.
\end{proof}

We obtain Theorem~\ref{NUthm}(a).

\begin{cor}\label{GQrigid}
$\mathcal G_\Q$ is conjugacy-rigid; in particular Theorem~\ref{NUthm}(a) holds.
\end{cor}

\begin{proof}
By Theorem~\ref{MainThm}(i) and Propositions~\ref{ExGQNeukirch} and \ref{PropInertialQ}
the group $\mathcal G_\Q$ is conjugacy-rigid. For Theorem~\ref{NUthm}(a) we may take
$K, K' \subseteq \overbar\Q$, as an embedding of $K$ in $\overbar\Q$ extends to an
isomorphism $\overbar K \simeq \overbar\Q$. Then
$\sigma\mathcal G_K\sigma^{-1} = \mathcal G_{K'}$ for some $\sigma \in \mathcal G_\Q$,
and taking fixed fields gives $K' = \sigma(K) \simeq K$.
\end{proof}

For Theorem~\ref{NUthm}(b) the arithmetic input is the following lemma.

\begin{lem}\label{LemRegQ}
Let $l$ be a prime and $F$ a finite Galois extension of $\Q$ with $\mu_l \subseteq F$
and $\Gamma = \Gal(F:\Q)$. There is a finite Galois extension $F_1$ of $\Q$ containing
$F$ such that $\Gal(F_1:F) \simeq \F_l[\Gamma]$ as $\Gamma$-modules, where $\Gamma$ acts
by conjugation. In particular $\mathcal G_\Q$ has enough regular layers.
\end{lem}

\begin{proof}
\textbf{Step 1. Construction of $F_1$.} Let $p$, $\mathfrak P$ and $a \in \mathcal O_F$ be
as in the proof of Lemma~\ref{LemQuad} for $K = N = F$. Then $v_{\mathfrak P}(a) = 1$. As
$p$ splits completely in $F$ the primes $\gamma\mathfrak P$ with $\gamma \in \Gamma$ are
distinct, so $v_{\mathfrak P}(\gamma(a)) = v_{\gamma^{-1}\mathfrak P}(a) = 0$ for
$\gamma \ne 1$. For
$\gamma \in \Gamma$ fix $\alpha_\gamma \in \overbar\Q$ with $\alpha_\gamma^l = \gamma(a)$
and put $F_1 = F(\alpha_\gamma : \gamma \in \Gamma)$. Let $\sigma \in \mathcal G_\Q$
restrict to $\gamma$. Then $\sigma(\alpha_\delta)^l = \gamma\delta(a)$, so
$\sigma(\alpha_\delta) \in \mu_l\,\alpha_{\gamma\delta}$ for $\delta \in \Gamma$. As
$\mu_l \subseteq F$ we get $\sigma(F_1) = F_1$, so $F_1$ is Galois over $\Q$. Also
$\rho \mapsto (\rho(\alpha_\delta)/\alpha_\delta)_{\delta \in \Gamma}$ is an injective
homomorphism from $\Gal(F_1:F)$ to $\mu_l^\Gamma$, as $\Gal(F_1:F)$ fixes $\mu_l$ and the $\alpha_\delta$ generate $F_1$ over $F$. So
$\Gal(F_1:F)$ is an $\F_l$-vector space of dimension at most $|\Gamma|$.

\medskip\noindent
\textbf{Step 2. An element $\tau \ne 1$ fixing $\alpha_\delta$ for $\delta \ne 1$.} Let $\Delta$ be the subgroup of $F^\times$ generated by $F^{\times l}$
and the $\delta(a)$ with $\delta \ne 1$. Every element of $\Delta$ has valuation at
$\mathfrak P$ divisible by $l$, so $a \notin \Delta$. By Kummer theory~\cite[Ch.~IV, Cor.~3.6]{NeukANT} distinct groups between $F^{\times l}$ and $F^\times$
have distinct fields of $l$-th roots over $F$. These fields are
$F_0 = F(\alpha_\delta : \delta \ne 1)$ for $\Delta$ and $F_1$ for the group generated by
$\Delta$ and $a$. So $F_0 \ne F_1$, and we take $\tau \ne 1$ in $\Gal(F_1:F_0)$.

\medskip\noindent
\textbf{Step 3. $\Gal(F_1:F)$ is the regular module $\F_l[\Gamma]$.} For $\gamma \in \Gamma$
let $\sigma \in \mathcal G_\Q$ restrict to $\gamma$ and put $\tau_\gamma =
\sigma\tau\sigma^{-1}$, which depends only on $\gamma$ as $\Gal(F_1:F)$ is abelian. For
$\delta \ne \gamma$ Step 1 gives $\alpha_\delta = c\,\sigma(\alpha_{\gamma^{-1}\delta})$
with $c \in \mu_l$, and $\tau$ fixes $\sigma^{-1}(c)\,\alpha_{\gamma^{-1}\delta}$ as $\gamma^{-1}\delta \ne 1$. So
$\tau_\gamma$ fixes $\alpha_\delta$, and as $\tau_\gamma \ne 1$ it moves $\alpha_\gamma$.
If a product of powers $\tau_\gamma^{n_\gamma}$ is trivial, then evaluating it at
$\alpha_\gamma$ shows that $l$ divides $n_\gamma$. By Step 1 the $|\Gamma|$ elements
$\tau_\gamma$ form a basis of $\Gal(F_1:F)$. As ${}^{\gamma'}\tau_\gamma =
\tau_{\gamma'\gamma}$ for $\gamma' \in \Gamma$, the map $e_\gamma \mapsto \tau_\gamma$ is an isomorphism
$\F_l[\Gamma] \to \Gal(F_1:F)$ of $\Gamma$-modules.

\medskip\noindent
\textbf{Step 4. Enough regular layers in $\mathcal G_\Q$.} The open normal subgroups of $\mathcal G_\Q$ are the groups $\mathcal G_E$ for finite
Galois extensions $E$ of $\Q$. For such $E$ the field $F = E(\mu_5)$ is finite Galois
over $\Q$ with $\mathcal G_F \subseteq \mathcal G_E$, and with $l = 5$ the layer
$\mathcal G_F/\mathcal G_{F_1}$ is regular. So $\mathcal G_\Q$ has enough regular layers.
\end{proof}

Theorem~\ref{MainThm}(ii), Proposition~\ref{ExGQNeukirch} and Lemma~\ref{LemRegQ}
give the following corollary. For number fields $K, K' \subseteq \overbar\Q$ it says
that every isomorphism $\mathcal G_K \to \mathcal G_{K'}$ is conjugation by some
$\sigma \in \mathcal G_\Q$, and then $\sigma(K) = K'$ by taking fixed fields. This is
Theorem~\ref{NUthm}(b), which is due to Uchida~\cite{Uchida1}.

\begin{cor}\label{CorStrongNU}
$\mathcal G_\Q$ is strongly conjugacy-rigid; in particular Theorem~\ref{NUthm}(b)
holds.
\end{cor}

With Theorem~\ref{CorFinite}(ii) in place of Theorem~\ref{MainThm} the zeta function
of one quadratic extension determines $K$, as follows.

\begin{cor}\label{CorZetaQuad}
Let $K$ be a number field and $a \in K$ as in Lemma~\ref{LemQuad}, and put
$L = K(\sqrt a)$. A number field $K'$ is isomorphic to $K$ if and only if some
quadratic extension $L'$ of $K'$ has $\zeta_{L'} = \zeta_L$.
\end{cor}

\begin{proof}
If $\tau : K \to K'$ is an isomorphism, then $L' = K'(\sqrt{\tau(a)})$ is isomorphic to
$L$, so $\zeta_{L'} = \zeta_L$. Conversely let $L'$ be a quadratic extension of $K'$
with $\zeta_{L'} = \zeta_L$. We may take $L, L' \subseteq \overbar\Q$. Let $F$ be the
Galois closure over $\Q$ of the compositum $LL'$ and $G = \Gal(F:\Q)$, and let $H$,
$U$, $H'$, $U'$ be the subgroups of $G$ fixing $K$, $L$, $K'$, $L'$. Then $U'$ is
normal of index $2$ in $H'$, and by Gassmann's criterion
$\Ind_U^G(1) \simeq \Ind_{U'}^G(1)$. The quadratic character of $\mathcal G_K$ with
kernel $\mathcal G_L$ is a Clifford character by the proof of Proposition~\ref{PropInertialQ}. As $\mathcal G_F \subseteq \mathcal G_L$, Lemma~\ref{LemCliffQuot}
shows that the quadratic character $\chi$ of $H$ with kernel $U$ satisfies
\eqref{Star}. Theorem~\ref{CorFinite}(ii) with $l = 2$ gives
$\sigma \in G$ with $H' = \sigma H\sigma^{-1}$. Taking fixed fields gives
$K' = \sigma(K) \simeq K$.
\end{proof}

\subsection{The solvable case}\label{SecSol}

For a number field $K \subseteq \overbar\Q$ let $K^{\mathrm{sol}} \subseteq \overbar\Q$
be the maximal prosolvable extension of $K$ and
$\mathcal G_K^{\mathrm{sol}} = \Gal(K^{\mathrm{sol}}:K)$ the maximal prosolvable
quotient of $\mathcal G_K$. The number fields
$K \subseteq \Q^{\mathrm{sol}}$ are those with solvable Galois closure $N$ over $\Q$.
For them $K^{\mathrm{sol}} = \Q^{\mathrm{sol}}$. Indeed $FK$ is solvable over $K$ for
every finite solvable Galois extension $F$ of $\Q$, and the Galois closure over $\Q$ of
a finite solvable Galois extension of $K$ is solvable over $N$ and hence over $\Q$. So
the group $\Gal(\Q^{\mathrm{sol}}:K) = \mathcal G_K^{\mathrm{sol}}$ is open of index
$[K:\Q]$ in $\mathcal G_\Q^{\mathrm{sol}}$, and every open subgroup is of this form.

As for $\mathcal G_\Q$, the first hypothesis is a theorem of Neukirch.

\begin{prop}[Gassmann-faithful]\label{PropSol}
$\mathcal G_\Q^{\mathrm{sol}}$ is Gassmann-faithful.
\end{prop}

\begin{proof}
Let $H = \Gal(\Q^{\mathrm{sol}}:K)$ and $H' = \Gal(\Q^{\mathrm{sol}}:K')$ be open
subgroups with $H \simeq H'$. As shown above they are $\mathcal G_K^{\mathrm{sol}}$ and
$\mathcal G_{K'}^{\mathrm{sol}}$. By a theorem of Neukirch~\cite[Satz~7]{NeukAnotherProof} the fields $K$ and $K'$ have the same Galois closure
$N$ over $\Q$, and $\Gal(N:K)$ and $\Gal(N:K')$ are Gassmann equivalent in
$\Gal(N:\Q)$. As $N$ lies in $\Q^{\mathrm{sol}}$, the permutation representations of $H$ and $H'$
are inflated from those of $\Gal(N:K)$ and $\Gal(N:K')$ in $\Gal(N:\Q)$, which are
isomorphic.
\end{proof}

Lemma~\ref{LemQuad} again gives the second hypothesis.

\begin{prop}[Enough Clifford characters]\label{PropSolMono}
$\mathcal G_\Q^{\mathrm{sol}}$ has enough Clifford characters.
\end{prop}

\begin{proof}
Let $\Gal(\Q^{\mathrm{sol}}:K)$ be an open subgroup and $a$ as in Lemma~\ref{LemQuad}. As shown at the start of this subsection
$K^{\mathrm{sol}} = \Q^{\mathrm{sol}}$, so $L = K(\sqrt a) \subseteq \Q^{\mathrm{sol}}$.
The Clifford character of $\mathcal G_K$ with kernel $\mathcal G_L$ from the proof of
Proposition~\ref{PropInertialQ} is trivial on the closed normal subgroup
$\mathcal G_{\Q^{\mathrm{sol}}} \subseteq \mathcal G_L$ of $\mathcal G_\Q$. So by
Lemma~\ref{LemCliffQuot} it defines a Clifford character of
$\Gal(\Q^{\mathrm{sol}}:K)$.
\end{proof}

We obtain the analogue of Corollary~\ref{GQrigid}.

\begin{cor}\label{SolvCor}
$\mathcal G_\Q^{\mathrm{sol}}$ is conjugacy-rigid; in particular, for number fields
$K, K' \subseteq \Q^{\mathrm{sol}}$ an isomorphism
$\Gal(\Q^{\mathrm{sol}}:K) \simeq \Gal(\Q^{\mathrm{sol}}:K')$ of topological groups implies $K \simeq K'$.
\end{cor}

\begin{proof}
As in the proof of Corollary~\ref{GQrigid}, using Theorem~\ref{MainThm}(i) and
Propositions~\ref{PropSol} and \ref{PropSolMono}.
\end{proof}

The strong form holds for the same fields. The open normal subgroups of
$\mathcal G_\Q^{\mathrm{sol}}$ are the groups $\Gal(\Q^{\mathrm{sol}}:E)$ for finite
Galois extensions $E$ of $\Q$ inside $\Q^{\mathrm{sol}}$. For such $E$ the field
$F = E(\mu_5)$ also lies in $\Q^{\mathrm{sol}}$, and the field $F_1$ of Lemma~\ref{LemRegQ} for $F$ and $l = 5$ is Galois over $\Q$ and abelian over $F$, so it lies
in $\Q^{\mathrm{sol}}$. Hence $\mathcal G_\Q^{\mathrm{sol}}$ has enough regular layers,
and Theorem~\ref{MainThm}(ii) with Proposition~\ref{PropSol} gives the following
corollary, which also follows from Uchida's theorem~\cite{Uchida1} with $\overbar\Q$ replaced by $\Q^{\mathrm{sol}}$, a replacement that works for every Galois extension of $\Q$ with no non-trivial abelian extension, as Iwasawa had noticed~\cite[p.~78]{NeukAst},~\cite[p.~617]{Uchida1}.

\begin{cor}\label{SolvStrongCor}
$\mathcal G_\Q^{\mathrm{sol}}$ is strongly conjugacy-rigid. In particular, for number
fields $K, K' \subseteq \Q^{\mathrm{sol}}$ every isomorphism
$\Gal(\Q^{\mathrm{sol}}:K) \to \Gal(\Q^{\mathrm{sol}}:K')$ of topological groups is conjugation by some
$\sigma \in \mathcal G_\Q^{\mathrm{sol}}$ with $\sigma(K) = K'$.
\end{cor}

\end{document}